\documentclass[11pt]{article}
\usepackage[margin=1in]{geometry}

\usepackage{amsmath} 
\usepackage{amssymb}  
\usepackage{amsthm}
\usepackage{color}

\usepackage{hyperref}
\usepackage{algorithm}

\theoremstyle{plain}  
\newtheorem{theorem}{Theorem}[section]

\newtheorem{proposition}[theorem]{Proposition}

\theoremstyle{definition}

\theoremstyle{remark} 

\title{\LARGE \bf Response to ``Counterexample to global convergence of DSOS and SDSOS hierarchies"}
 \author{Amir Ali Ahmadi\thanks{Department of Operations Research and Financial Engineering, Princeton University.
       ({a\_a\_a@princeton.edu}, \url{http://aaa.princeton.edu/})} and Anirudha Majumdar  \thanks{Department of Mechanical and Aerospace Engineering, Princeton University. ({ani.majumdar@princeton.edu}, \url{https://irom-lab.princeton.edu/majumdar}).} }
\begin{document}
\date{}
\maketitle


\begin{abstract}


In a recent note~\cite{cedric}, the author provides a counterexample to the global convergence of what his work refers to as ``the DSOS and SDSOS hierarchies'' for polynomial optimization problems (POPs) and purports that this refutes claims in our extended abstract~\cite{dsos_ciss} and slides in~\cite{dsos_slides}. The goal of this paper is to clarify that neither~\cite{dsos_ciss}, nor~\cite{dsos_slides}, and certainly not our full paper~\cite{dsos}, ever defined DSOS or SDSOS hierarchies as it is done in~\cite{cedric}. It goes without saying that no claims about convergence properties of the hierarchies in~\cite{cedric} were ever made as a consequence. What was stated in~\cite{dsos_ciss,dsos_slides} was completely different: we stated that \emph{there exist} hierarchies based on DSOS and SDSOS optimization that converge. This is indeed true as we discuss in this response. We also emphasize that we were well aware that some (S)DSOS hierarchies do not converge even if their natural SOS counterparts do. This is readily implied by an example in our prior work~\cite{dsos}, which makes the counterexample in~\cite{cedric} superfluous. Finally, we provide concrete counterarguments to claims made in \cite{cedric} that aim to challenge the scalability improvements obtained by DSOS and SDSOS optimization as compared to sum of squares (SOS) optimization. 

\end{abstract}



\section{Introduction}

In \cite{cedric}, the author considers our recently proposed (S)DSOS optimization framework~\cite{dsos}. Based on it, he constructs two particular hierarchies, which \cite{cedric} calls ``the DSOS and SDSOS hierarchies'', and shows that they do not always converge to the global optimal value of polynomial optimization problems (POPs) via an example. The author then makes two primary claims: (i) that this example refutes statements made in our extended abstract \cite{dsos_ciss} and slides \cite{dsos_slides}, and (ii) that DSOS and SDSOS optimization do not provide tractable alternatives to SOS optimization. This response shows that (i) is not true and the main technical premise on which (ii) is based is also not true.\footnote{We take this opportunity to point out that the reference for our work on DSOS and SDSOS optimization should be our full paper on the topic~\cite{dsos}. In both our slides \cite{dsos_slides}, and our extended abstract~\cite{dsos_ciss}, details have intentionally been omitted (see e.g. the disclaimer on the first page of~\cite{dsos_ciss}).}

We refer the reader to \cite{dsos} for a formal introduction to DSOS and SDSOS programming.\footnote{Our notation here is consistent with \cite{dsos}; in particular, we use upper case to refer to (S)DSOS cones and hierarchies, while lower case is reserved for (s)dsos polynomials.} Briefly, dsos (resp. sdsos) polynomials are sos polynomials whose Gram matrices are diagonally dominant (resp. scaled diagonally dominant). Optimizing a linear function over these cones subject to affine constraints can be performed using linear and second-order cone programming respectively. The resulting DSOS and SDSOS programs have shown improvements in scalability over SOS programming for a range of different problems \cite[Section 4]{dsos}.

\section{Claim 1: Non-convergence of DSOS and SDSOS hierarchies on POPs}
In reference to our work, the abstract of \cite{cedric} states that the result presented in \cite{cedric} ``refutes the claim in the literature according to which the DSOS and SDSOS hierarchies can solve any polynomial optimization problem to arbitrary accuracy".

\subsection{Rebuttal to Claim 1}

The ``claim in the literature" that \cite{cedric} purports to refute does not appear anywhere in our main manuscript~\cite{dsos}, which was made publicly available in 2017. It does not even feature (contrary to what~\cite{cedric} claims) in our extended abstract~\cite{dsos_ciss}  from 2014 or presentation slides~\cite{dsos_slides} from 2013.
%
In particular, our prior work did not present ``the DSOS and SDSOS hierarchies" that \cite{cedric} considers to be convergent hierarchies for the POP. 
The results in \cite{cedric} therefore do not provide counterexamples to any results proven in \cite{dsos, dsos_ciss} (or, to our knowledge, any other claims made in the literature).

What we do claim in~\cite{dsos_ciss,dsos_slides} is that there \emph{exist} converging hierarchies for POPs based on DSOS and SDSOS optimization. This statement is certainly true; in fact there are many weaker converging hierarchies that as an immediate corollary imply existence of a (S)DSOS-based converging hierarchy (see, e.g., those by Lasserre~\cite{Lasserre_book},~\cite[Chap. 9]{Lasserre_book2}, Pe\~na, Vera, and Zuluaga~\cite{pena_LP},~\cite[Theorem 2]{kuang_alternative}, and Ahmadi and Hall~\cite{pop_hierarchy}). In the first two references, the LP hierarchy presented only requires a search over nonnegative scalars. If one replaces scalar optimization variables by dsos polynomials, one trivially obtains new LP-based globally convergent hierarchies (since the modified hierarchy automatically inherits the convergence properties of the original LP hierarchy). In the latter reference, the sos polynomials that feature in the hierarchy all have a \emph{diagonal} Gram matrix. These polynomials are also a very special case of dsos polynomials. The assumptions that these hierarchies make on the POP (particularly in the case of~\cite{pop_hierarchy} and~\cite{Lasserre_book2}) are very similar to those made for SOS-based converging hierarchies; see, e.g., the discussion in~\cite[Section 2.4.3]{Lasserre_book2}. In general, we expect that there are many more possible ways of constructing hierarchies based on (S)DSOS optimization. We further note that \cite[Section 3.2]{dsos} was already providing a convergent hierarchy for the important special case of copositive programs~\cite{burer_copositive}.

%
%


Finally, we note that while globally convergent hierarchies are valuable from a theoretical perspective, they are generally of limited utility to practitioners without concrete and practical bounds on convergence rates. This is a challenge for the usual hierarchies based on SOS optimization, and it is also a challenge faced by convergent hierarchies based on (S)DSOS optimization. Due to this reason, we did not focus on explicitly presenting convergent hierarchies for POP in \cite{dsos}, which is meant to be a broadly-accessible application-oriented paper. Rather, in \cite{dsos}, we chose to focus on demonstrating the scalability of (S)DSOS programming using large-scale numerical examples drawn from a diverse range of application areas (\cite[Section 4]{dsos}), and algorithms for iteratively improving the quality of solutions obtained from (S)DSOS programs (\cite[Section 5]{dsos}). In a different paper~\cite{pop_hierarchy}, one can find a systematic approach for constructing converging hierarchies for POPs based on minimal requirements. In fact, the methodology presented there produces a converging hierarchy that does not even use optimization.

\subsection{An example from our prior work makes the main claim in~\cite{cedric} superfluous}

We point out in this subsection that our paper~\cite{dsos}, which is referenced in~\cite{cedric}, already gave an example of a very simple polynomial optimization problem (see problem (\ref{eq:simple.pop}) below) for which the first level of a well-known SOS hierarchy is exact but for which no level of the analogous DSOS or the SDSOS hierarchies are exact. This is an immediate corollary of~\cite[Proposition 14]{dsos} in our paper. It shows that we were cognizant of the fact that one cannot take any SOS hierarchy, replace the sos constraints with dsos or sdsos constraints (or even $r$-dsos and $r$-sdsos constraints), and expect convergence. As a consequence, we would not have made such a claim, and we certainly never stated that the Lasserre hiearchy with sos constraints replaced by dsos/sdsos constraints converges. Our claim was  simply that one can construct hierarchies with dsos/sdsos constraints that converge. 

\vspace{5mm}

Recall that a polynomial $f$ is said to be $r$-sos (resp. $r$-dsos, $r$-sdsos) if $f(x)\cdot(\sum x_i^2)^r$ is sos (resp. dsos, sdsos). Consider the problem of minimizing a form $p$ of degree $2d$ on the sphere:
\begin{equation}\label{eq:p.on.sphere}
\begin{aligned}
p^*\mathop{\mathrel:}=&\min_{x \in \mathbb{R}^n} && p(x)\\
&\text{s.t. } && \sum_{i=1}^n x_i^2=1.\\
\end{aligned}
\end{equation}

The most standard SOS hierarchy for this problem (indexed by $r\in\mathbb{N}$) reads~\cite{PhD:Parrilo}:

\begin{equation}\nonumber
\begin{aligned}
l_{sos,r}\mathop{\mathrel:}=&\max_{\gamma \in \mathbb{R}} && \gamma\\
&\text{s.t. } &&  p(x)-\gamma(\sum_{i=1}^n x_i^2)^d \mbox{\quad is $r$-sos.}\\
\end{aligned}
\end{equation}

These are semidefinite programs that produce a sequence of lower bounds $\{l_{sos,r}\}$ on $p^*$ with $l_{sos,r}\rightarrow p^*$ as $r\rightarrow\infty$. This hierarchy is shown to be equivalent to the Lasserre hierarchy applied to problem (\ref{eq:p.on.sphere}) in~\cite[Proposition 2]{pablo_etienne_monique}.

One can now define an analogous DSOS/SDSOS hierarchy:

\begin{equation}\nonumber
\begin{aligned}
l_{dsos,r} \mbox{\ (resp. $l_{sdsos,r}$)}\mathop{\mathrel:}=&\max_{\gamma \in \mathbb{R}} && \gamma\\ 
&\text{s.t. } &&  p(x)-\gamma(\sum_{i=1}^n x_i^2)^d \mbox{\quad is $r$-dsos (resp. $r$-sdsos).}\\
\end{aligned}
\end{equation}
We clearly have $l_{dsos,r}\leq l_{sdsos,r} \leq l_{sos,r}$ for all $r$. Consider now the following optimization problem:
\begin{equation}\label{eq:simple.pop}
\begin{aligned}
&\min_{x \in \mathbb{R}^3} && (x_1+x_2+x_3)^2\\
&\text{s.t. } && x_1^2+x_2^2+x_3^2=1.\\
\end{aligned}
\end{equation}

The optimal value of this problem is zero, which is achieved by the first-order sos relaxation; i.e., $l_{sos,0}=0$.  (This is because all nonnegative quadratic forms are sums of squares.) However, one has $$l_{dsos,r}\leq l_{sdsos,r}\leq -1$$ for all $r$. This is an immediate consequence of the following proposition, which we have already proven in~\cite{dsos}.

\begin{proposition}[Proposition 14 in \cite{dsos}]
For any $0<a<1$, the quadratic form
\begin{equation}
 f(x_1,x_2,x_3)=(x_1+x_2+x_3)^2+a(x_1^2+x_2^2+x_3^2)
\end{equation}
is positive definite but not r-sdsos for any $r$.
\end{proposition}

We remark that in contrast to our example in (\ref{eq:simple.pop}), the example given in~\cite{cedric} is in two variables and constitutes a convex problem, but these considerations are secondary.

\section{Claim 2: (S)DSOS is not necessarily more tractable than SOS}
The abstract of \cite{cedric} states: ``We further observe that the dual to the SDSOS hierarchy is the moment hierarchy where every positive semidefinite constraint is relaxed to all necessary second-order conic constraints. As a result, the number of second-order conic constraints grows exponentially as the order of the SDSOS hierarchy increases. Together with the counterexample, this suggests that DSOS and SDSOS are not necessarily more tractable alternatives to sum-of-squares."

\subsection{Rebuttal to Claim 2}
As observed already in earlier work \cite[Section 3.2]{permenter_parrilo_facial_reduction},~\cite[Section 3.3]{basis_pursuit_Ahmadi_Hall}, the dual to SDSOS programs can indeed be obtained by relaxing semidefinite constraints to necessary second-order conic constraints. However, the number of SOCP constraints grows \emph{polynomially} as the order of the $r$-sdsos hierarchy increases. Hence the claim above about the exponential increase is false. Indeed, consider a polynomial in $n$ variables and degree $2d$. When this polynomial is required to be $r$-sdsos, the scaled diagonal dominance constraint needs to be imposed on a symmetric matrix of size $N\times N$, where $N\mathrel{\mathop:}=\binom{n+d+r}{n}$. This leads to $O(N^2)=O(n+d+r)^{2n}$ SOCP constraints. Hence, the dependence on $r$ is indeed polynomial. (Here, $n$ is constant as the number of variables to any polynomial optimization problem that goes through the hierarchy is clearly fixed.)

We shall emphasize that the size of the matrices that SDSOS optimization deals with is exactly the same as that of the SOS approach. The difference is that an expensive semidefinite constraint of size $N\times N$ is replaced with $\binom{N}{2}$ cheap SOCP constraints. We have shown with numerous examples (see~\cite[Section 4]{dsos}) that for a range of problem sizes of practical value, this can lead to significant improvements in scalability. We also refer the reader to Theorem 10 and Theorem 12 of~\cite{dsos}, where we show that from a theoretical standpoint, polynomial-time solvability of $r$-(S)DSOS programs is identical to that of $r$-SOS programs.

%
%
%


\section{Discussion}
From an application viewpoint, we believe that the (S)DSOS optimization approach can provide practitioners in diverse application domains with powerful and tractable alternatives to SOS optimization. As we demonstrate in \cite{dsos, Majumdar14a} with numerical examples from a wide range of domains (including copositive and combinatorial optimization, machine learning and statistics, control theory, and robotics), the (S)DSOS approach can provide significant gains in computation times (as much as orders of magnitude in certain cases) when compared to the SOS approach. For example, in \cite{Majumdar14a}, we used SDSOS programming to design stabilizing feedback controllers for a $30$ dimensional state space and $14$ dimensional control input space model of a humanoid robot\footnote{A video is available online: \href{http://youtu.be/lmAT556Ar5c}{http://youtu.be/lmAT556Ar5c}}. The scale of this problem is well beyond what SOS programming can currently handle. The computational gains in the problems considered in \cite{dsos, Majumdar14a} are obtained without exploiting any particular structure (e.g., sparsity or symmetry) in the problems. As noted in \cite[Section 2]{dsos}, the (S)DSOS approach could potentially be combined with approaches that exploit the structure of the problem at hand in order to obtain even more significant computational gains. We further note that while there is conservatism inherent in the (S)DSOS approach as compared to SOS optimization, the examples in \cite{dsos} demonstrate that this conservatism can be small and mitigated by the techniques presented in~\cite[Section 5]{dsos}. 

We conclude by noting that the (S)DSOS framework is not meant as a general replacement for hierarchies based on SOS programming (e.g., the Lasserre/Parrilo hierarchies). Rather, the (S)DSOS approach provides a way to tackle problems that are beyond the reach of SOS programming due to computational limitations. Indeed, we believe that the true power of the (S)DSOS framework comes from the fact that it can provide solutions in situations where even the first level of the SOS relaxation hierarchy is simply too expensive to solve (see \cite{dsos, Majumdar14a} for numerous examples of this kind).



\bibliographystyle{abbrv}
\bibliography{pablo_amirali}





\end{document}